\documentclass{elsart}
\usepackage{amsfonts,amsmath}
\usepackage{amssymb,latexsym,amsmath}

\begin{document}

\begin{frontmatter}
\title{Krall--type Orthogonal Polynomials \\in Several Variables\thanksref{td}}
\subtitle{\small Dedicated to J. S. Dehesa on the occasion of his 60th birthday}

\author{Lidia Fern\'{a}ndez, Teresa E. P\'erez, Miguel A. Pi\~{n}ar}
\address{Departamento de Matem\'atica Aplicada,
 and Instituto Carlos I de F\'{\i}sica Te\'orica y Computacional,
Universidad de Granada (Spain)} \author{Yuan Xu}
\address{Department of Mathematics,
  University of Oregon (USA)}
\thanks[td]{Partially supported by Ministerio de Ciencia y Tecnolog\'{\i}a
        (MCYT) of Spain and by the European Regional Development Fund
        (ERDF) through the grant MTM 2005--08648--C02--02, and Junta de
        Andaluc\'{\i}a, Grupo de Investigaci\'on FQM 0229.  The work of the
         fourth author is supported in part by NSF Grant DMS--0604056.}

\begin{abstract}
For a bilinear form obtained by adding a Dirac mass to a positive definite
moment functional in several variables, explicit formulas of orthogonal
polynomials are derived from the orthogonal polynomials associated with
the moment functional. Explicit formula for the reproducing kernel is also
derived and used to establish certain inequalities for classical orthogonal
polynomials.

\end{abstract}

\end{frontmatter}

\section{Introduction}

Let $\Pi^d$ denote the space of polynomials in $d$--variables, and
let $u$ be a moment functional, denoted by $\langle u, p \rangle$
for $p \in \Pi^d$, for which orthogonal polynomials exist. We
define a new functional $v$ by adding a Dirac mass to $u$,
\begin{equation}\label{v-u}
 \langle v, p\rangle = \langle u, p\rangle + \lambda \, p(c), \qquad c \in \mathbb{R}^d,
     \quad \lambda \in  \mathbb{R}, \quad p \in \Pi^d,
\end{equation}
and study orthogonal polynomials with respect to the functional $v$.

In the case of one variable, this problem first arose from the work of A. M. Krall
(\cite{K1}) when he studied the orthogonal polynomials that are eigenfunctions of
a fourth order differential operator considered by H. L. Krall (\cite{K2,K3}), and
showed that the polynomials are orthogonal with respect to a measure that is
obtained from a continuous measure on an interval by adding masses at the end
points of the interval.  In \cite{Koor}, Koornwinder studied the case that the
measure is a Jacobi weight function together with additional mass points at
$1$ and $-1$; he constructed explicit
orthogonal polynomials and studied their properties. Uvarov (\cite{U}) considered
the problem of orthogonal polynomials with respect to a measure obtained by
adding a finite discrete part to another measure; his main result expresses the
polynomials orthogonal with respect to the new measure in terms of the
polynomials orthogonal with respect to the old one. More generally, one can
consider perturbations of quasi definite linear functionals via the addition of
Dirac delta functionals, orthogonal polynomials in such a general setting has
been studied extensively in recent years (see, for instance \cite{MM}, and the
references therein).

The purpose of the present paper is to study this problem in several variables.
After a brief section on notations and preliminaries in the next section, we
state and prove our main results in Sections 3. The result gives a necessary and
sufficient condition for the existence of orthogonal polynomials with respect to the
linear functional $v$ defined by \eqref{v-u}, and expresses orthogonal polynomials
with respect to $v$ in terms of the orthogonal polynomials with respect to $u$. Furthermore, we can also express the reproducing kernel of polynomials with
respect to $v$ in terms of the kernel with respect to $u$. Our formula on the
reproducing kernel implies an inequality on the orthogonal polynomials with
respect to $u$. Even in the case of one variable, it leads to new inequalities on
classical orthogonal polynomials, which are stated in Section 4. Finally, in
Section 5, we consider the example of orthogonal polynomials on the unit ball.

\section{ Orthogonal polynomials in several variables}

In this section we recall necessary notations and definitions about orthogonal
polynomials of several variables, following \cite{DX}.

Throughout this paper, we will use the usual multi--index
notation. For $\alpha=(\alpha_1,\dots,\alpha_d) \in\mathbb{N}_0^d$
and $x=(x_1,\dots,x_d) \in\mathbb{R}^d$, we write $x^{\alpha}=
x_1^{\alpha_1}\cdots x_d^{\alpha_d}$. The integer
$|\alpha|=\alpha_1+ \dots+\alpha_d$ is called the \emph{total
degree} of $x^{\alpha}$. The linear combinations of $x^\alpha$,
$|\alpha| =n$, is a homogeneous polynomial of degree $n$.  We
denote by $\mathcal{P}_n^d$ the space of homogeneous polynomials
of degree $n$ in $d$ variables, by $\Pi_n^d$ the space of
polynomials of total degree not  greater than $n$. It is well
known that
$$
\dim \Pi_n^d = \binom{n+d}{n} \quad \hbox{\rm and}\quad \dim
\mathcal{P}_n^d = \binom{n+d-1}{n}:=r_n^d.
$$
Let $\{\mu_{\alpha}\}_{\alpha\in \mathbb{N}_0^d}$ be a multi-sequence of
real numbers, and let $u$ be a real valued functional defined on
${\mathcal P}$ by means of
$$
\langle u,x^{\alpha}\rangle = u (x^{\alpha}) = \mu_{\alpha},
$$
and extended by linearity. Then, $u$ is called the \emph{moment
functional} determined by $\{\mu_{\alpha}\}_{\alpha\in
\mathbb{N}_0^d}$. If $\langle u, p^2\rangle > 0$, $\forall p\in
\Pi_n^d$ and $p\neq 0$, then the moment $u$ is called
\emph{positive definite} and it induces an inner product
accordingly by $(p, q) := \langle u, p\, q\rangle$, $\forall p,
q\in \Pi^d. $ A typical example of a positive moment functional is
an integral with respect to a positive measure $d\mu$ with all
moments finite, $\langle u, p \rangle = \int_{\mathbb{R}^d} p(x)
d\mu$.

A polynomial $P \in \Pi_n^d$ is called an orthogonal polynomial
with respect to $u$ if $\langle u, P Q \rangle =0$ for all $Q \in
\Pi_{n-1}^d$. Let $V_n^d$ denote the space of orthogonal polynomials
with respect to the $u$.  We are interested in the case when $u$ admits
a basis of orthogonal polynomials; that is, $\dim V_n^d = r_n^d$ for
all $n$. This happens whenever $u$ is positive definite. In general, we
call a moment functional {\it quasi definite} if it admits a basis of
orthogonal polynomials.

For the study of orthogonal polynomials it is often convenient to  adopt
a vector notation (\cite{Ko1,X93}). Let $\{P_\alpha^n\}_{|\alpha|=n}$
denote a basis of $V_n^d$. Let the elements of $\{\alpha \in \mathbb{N}^d:
|\alpha| = n\}$ be ordered by $\alpha_1, \alpha_2, \ldots,  \alpha_{r_n^d}$
according to a fixed monomial order, say the lexicographical order. We
then write the basis of $V_n^d$ as a column vector
$$
\mathbb{P}_n=(P_{\alpha}^n)_{|\alpha|=n}=
 (P_{\alpha_1}, P_{\alpha_2}, \ldots,  P_{\alpha_{r_n^d}})^T.
$$
Using this notation, the orthogonality of $P_\alpha^n$ can be expressed as
\begin{equation*} 
  \langle u,{\mathbb{P}}_{n}{\mathbb{P}}_{m}^T\rangle = H_n \delta_{m,n},
\end{equation*}
where $H_n$ is a matrix of size $r_n^d\times r_n^d$. That $u$ is
quasi definite is equivalent to that $H_n$ is invertible for all $n \ge 0$.
If $u$ is positive definite, then we can choose a basis so that $H_n$ is
the identity matrix for all $n$;  in other words, we can choose the basis
to be {\it orthonormal}.

Let $u$ be a quasi definite moment functional and let
$\{\mathbb{P}_n\}_{n\ge 0}$ be a sequence of orthogonal polynomials
with respect to $u$. The reproducing kernel of $V_n^d$ and $\Pi_n^d$
are denoted by $\mathbf{P}_n(\cdot, \cdot)$ and $\mathbf{K}_n(\cdot, \cdot)$,
respectively, and are given by
$$
  \mathbf{P}_k(x,y) = \mathbb{P}^T_k(x)\, H^{-1}_k\, \mathbb{P}_k(y) \quad
\hbox{and} \quad
\mathbf{K}_n(x,y) = \sum_{k=0}^n \mathbf{P}_k(x,y).
$$
These kernels satisfy the usual reproducing property; for example,
$$
   p(x) = \langle u, \mathbf{K}_n(x, \cdot)\, p(\cdot) \rangle
      = \langle u, p(\cdot)\, \mathbf{K}_n(\cdot,
    x) \rangle, \quad p\in \Pi_n^d,
$$
which shows, in particular, that these functions are independent
of a particular choice of the bases. The kernel
$\mathbf{K}_n(x,y)$ satisfies an analog of the
Christoffel--Darboux formula, and plays an important role in the
study of orthogonal Fourier expansions. For this and further
properties of orthogonal polynomials of several variables, see
\cite{DX}.

\section{Krall--type orthogonal polynomials in several variables}

Let $u$ be a quasi definite moment functional defined on $\Pi^d$.
We define a new moment functional $v$ as the perturbation of $u$ given by
$$
\langle v, p\rangle = \langle u, p\rangle + \lambda \,
p(c), \quad \forall p(x)\in \Pi^d,
$$
where $\lambda$ is a non zero real number and $c\in \mathbb{R}^d$ is a
given point. Our first result gives a necessary and sufficient condition for
$v$ to be quasi definite.

\begin{thm} \label{main-thm}
The moment functional $v$ is quasi definite if and only if
$$
\lambda_n := 1\, + \, \lambda \, \mathbf{K}_{n}(c,c) \neq 0, \quad n\ge 0.
$$
Furthermore, when $v$ is quasi definite, a sequence of orthogonal
polynomials, $\{\mathbb{Q}_n\}_{n\ge0}$, with respect to $v$ is
given by
\begin{equation}\label{ex-expl}
\mathbb{Q}_n(x) = \mathbb{P}_n(x)\, -
\,\frac{\lambda}{\lambda_{n-1}}\,\mathbf{K}_{n-1}(c,x)\,\mathbb{P}_n(c)
 \quad n\ge0,
\end{equation}
where $\{\mathbb{P}_n\}_{n\ge0}$ denote orthogonal polynomials
with respect to $u$ and $\mathbf{K}_{-1} (\cdot,\cdot) :=0$.
\end{thm}

\begin{pf}
First we assume that $v$ is quasi definite and
$\{\mathbb{Q}_n\}_{n\ge0}$ is a sequence of orthogonal polynomials
with respect to $v$. Since $u$ is quasi definite, the leading
coefficient of $\mathbb{P}_n$ is an invertible matrix. Hence, by
multiplying $\mathbb{Q}_n$ by an invertible matrix, if necessary,
we can assume that $\mathbb{Q}_n-\mathbb{P}_n \in
\mathcal{P}_{n-1}$ for $n\ge 0$. This shows, in particular, that
$\mathbb{Q}_0 = \mathbb{P}_0$. Furthermore, since
$\{\mathbb{P}_n\}_{n \ge 0}$ is a basis of $\Pi^d$, for each $n\ge
1$ we can express $\mathbb{Q}_n$ in terms of $\mathbb{P}_n$. Thus,
there exist constant matrices $M_i^n$ of size $r_n^d\times r_i^d$
such that
$$
\mathbb{Q}_n(x) = \mathbb{P}_n(x) +
  \sum_{i=0}^{n-1} M_i^n\, \mathbb{P}_i(x),
$$
where, by the orthogonality of $\mathbb{P}_n$ and the definition of $v$,
$$
M_i^n = \langle u, \mathbb{Q}_n\, \mathbb{P}^T_i\rangle \, H^{-1}_i= \left[\langle v,
\mathbb{Q}_n\,\mathbb{P}^T_i \rangle - \lambda\, \mathbb{Q}_n(c)\,
\mathbb{P}^T_i(c)\right] H^{-1}_i,
$$
for $0\le i\le n-1$. Since $\langle v,\mathbb{Q}_n\,\mathbb{P}^T_i \rangle =0$
for $i \le n-1$, we conclude then
$$
\mathbb{Q}_n(x) = \mathbb{P}_n(x) - \lambda\, \mathbb{Q}_n(c)\,
\sum_{i=0}^{n-1} \mathbb{P}^T_i(c) \, H^{-1}_i \, \mathbb{P}_i(x)
= \mathbb{P}_n(x) - \lambda\, \mathbb{Q}_n(c)\,
\mathbf{K}_{n-1}(c,x).
$$
Evaluating the above expression at $x=c$, we obtain
\begin{equation}\label{qn(c)}
\mathbb{Q}_n(c) \left[1 + \lambda\, \mathbf{K}_{n-1}(c,c)\right] =
  \mathbb{P}_n(c).
\end{equation}
Recall that  $\lambda_k = 1 + \lambda\, \mathbf{K}_k(c,c)$. If $\lambda_{n-1}
 = 0$ for some value of $n\ge 1$, then $\mathbb{P}_n(c)=0$ by \eqref{qn(c)}
 and, furthermore, $\lambda_n=\lambda_{n-1} + \lambda\,\mathbb{P}_n^T(c)\,
  H^{-1}_n\, \mathbb{P}_n(c)=0$. Thus, we conclude that
$\mathbb{P}_{n+1}(c)=0$ which, however, contradict to the fact
that $\mathbb{P}_{n}$ and $\mathbb{P}_{n+1}$ cannot have a common
zero (\cite{DX}, p. 113). This shows that $ \lambda_n\neq 0$,
$\forall n\ge0$, and also that (\ref{ex-expl}) holds whenever
$\lambda_n \ne 0$.

Conversely, if all $\lambda_n$ are non--zero, then we can define
the polynomials $\mathbb{Q}_n$  by \eqref{ex-expl} and the above
proof shows then that $\mathbb{Q}_n$ is orthogonal with respect to
$v$. Since $\mathbb{Q}_n$ and $\mathbb{P}_n$ has the same leading
coefficient matrix, it is evident that $\{\mathbb{Q}_n \}_{n \ge
0}$ contains a basis of $\Pi^d$. \qed\end{pf}

\begin{rem}
If $u$ is a quasi definite moment functional, then $v$ is quasi definite
except when $\lambda$ belongs to an infinite discrete set of values. If
$u$ is positive definite and we choose $\lambda > 0$, then it is easy to
see that $v$ is also positive definite (see \eqref{inversa} below).
\end{rem}

Let $\mathbb{Q}_n$ be as in the theorem, we define $\widetilde H_k
= \langle v, \mathbb{Q}_k \, \mathbb{Q}_k^T\rangle$ and  denote by
\begin{equation*}
\mathbf{\widetilde P}_k(x,y):= \mathbb{Q}^T_k(x)\, \widetilde
H^{-1}_k\, \mathbb{Q}_k(y) \quad \hbox{and} \quad
\mathbf{\widetilde K}_n(x,y):= \sum_{k=0}^n \mathbf{\widetilde
P}_k(x,y),
\end{equation*}
the reproducing kernels associated with the linear functional $v$.  We will
derive an explicit formula of these kernels when $u$ is positive definite.

Whenever $u$ is positive definite, we choose $\{\mathbb{P}_n\}_{n
\ge 0}$ in the Theorem \ref{main-thm} as a sequence of orthonormal
basis with respect to $u$, so that $H_n  = \langle u,
\mathbb{P}_n\mathbb{P}_n^T\rangle =I_{r_n}$, the identity matrix
of size $r_n$.

\begin{prop}
Let $u$ be a positive definite linear functional. Then, for $k \ge 0$,
\begin{equation}\label{inversa}
 \widetilde H_k = I_{r_k} + \frac{\lambda}{\lambda_{k-1}}
    \mathbb{P}_k(c) \, \mathbb{P}_k (c)^T \quad\hbox{and} \quad
\widetilde H^{-1}_k = I_{r_k} - \frac{\lambda}{\lambda_{k}}
\mathbb{P}_k (c) \,
    \mathbb{P}_k(c)^T.
\end{equation}
\end{prop}

\begin{pf}
Since we choose $\mathbb{P}_n$ so that $H_n = I_{r_n}$, it follows
from (\ref{qn(c)}) that
\begin{eqnarray*}
\widetilde H_k = \langle v, \mathbb{Q}_k \, \mathbb{Q}_k^T\rangle =
\langle v, \mathbb{Q}_k \, \mathbb{P}_k^T\rangle &=& \langle u,
\mathbb{Q}_k \, \mathbb{P}_k^T\rangle + \lambda \,\mathbb{Q}_k(c)\,
\mathbb{P}_k(c)^T \\
&=& I_{r_k} + \frac{\lambda}{\lambda_{k-1}} \,\mathbb{P}_k(c) \,
\mathbb{P}_k(c)^T.
\end{eqnarray*}
Assuming that the inverse is of the form $\widetilde H^{-1}_k =
I_{r_k} - \delta \mathbb{P}_k (c) \, \mathbb{P}_k(c)^T$, and using
the fact that $\mathbb{P}_k (c)^T \mathbb{P}_k(c) =
\mathbf{P}_n(c,c)$, a quick computation shows that $\widetilde H_k
\widetilde H_k^{-1} = I_{r_k}$ is equivalent to $\lambda -
\lambda_{k-1} \delta - \lambda \delta \mathbf{P}_n(c,c) = 0$.
Using the fact that $\mathbf{P}_n(c,c)= \mathbf{K}_n(c,c) -
\mathbf{K}_{n-1}(c,c)$, it is easy to see that $\delta = \lambda /
\lambda_{k}$.
\end{pf}

\begin{thm}
Let $u$ be a positive definite linear functional. Then, for $k\ge 0$,
\begin{equation} \label{Pk}
   \mathbf{ \widetilde P}_k(x,y) = \mathbf{P}_k(x,y) - \frac{\lambda}{\lambda_k}\,
\mathbf{K}_k(x,c)\,\mathbf{K}_k(c,y) +
\frac{\lambda}{\lambda_{k-1}}\,\mathbf{K}_{k-1}(x,c)\,
\mathbf{K}_{k-1}(c,y).
\end{equation}
Furthermore, for $n \ge 0$,
\begin{equation}\label{kernel}
\mathbf{\widetilde K}_n(x,y) =
\mathbf{K}_n(x,y) -
\frac{\lambda}{\lambda_n}\,\mathbf{K}_n(x,c)\,
\mathbf{K}_n(c,y).
\end{equation}
\end{thm}

\begin{pf}
Let $\gamma_k: = \lambda/ \lambda_k$. By \eqref{ex-expl} and \eqref{inversa},
it follows readily that
\begin{align*}
 \mathbf{\widetilde P}_k(x,y) = \, & \mathbf{P}_k(x,y) -  \gamma_k\left[
  \mathbf{P}_k(x,c) \mathbf{K}_{k-1}(x,c)+  \mathbf{P}_k(y,c)\mathbf{K}_{k-1}(c,x) \right. \\
    & \left. \qquad  +  \mathbf{P}_k(x,c) \mathbf{P}_k(y,c)\right]
       + \gamma_k \gamma_{k-1} \mathbf{P}_k(c,c) \mathbf{K}_{k-1}(x,c)\mathbf{K}_{k-1}(y,c)\\
   = & \, \mathbf{P}_k(x,y) -  \gamma_k \mathbf{K}_k(x,c) \mathbf{K}_{k}(y,c) \\
      & \qquad + \gamma_k [1+ \gamma_{k-1} \mathbf{P}_k(c,c)] \mathbf{K}_{k-1}(x,c)
           \mathbf{K}_{k-1}(y,c),
\end{align*}
as the first square bracket is equal to $\mathbf{K}_{k}(x,c)\mathbf{K}_{k}(y,c)-
\mathbf{K}_{k-1}(x,c)\mathbf{K}_{k-1}(y,c)$. Since the definition of $\gamma_k$ leads
readily to $1+ \gamma_{k-1} \mathbf{P}_k(c,c) = \gamma_{k-1}$, this proves \eqref{Pk}.
Summing over \eqref{Pk} for $k = 0, 1, \ldots, n$ proves \eqref{kernel}.
\qed \end{pf}
We note that, in the case of one variable, the formula \eqref{kernel} has appeared in
\cite{GPRV}, whereas the formula \eqref{Pk} appears to be new even in one variable.

\section{An application of formula \eqref{Pk}}

If $u$ is positive definite and $\lambda > 0$, then $v$ is also
positive definite. As a result, $\widetilde H_n$, hence
$\widetilde H^{-1}_n$ are positive definite matrices by
\eqref{inversa} for $n\ge0$. In particular, $\mathbf{\widetilde
P}_n(x,x)$ is nonnegative. In fact, if $d \ge 2$ and the linear
functional $u$ is centrally symmetric (see \cite{DX} for
definition), then $\mathbf{\widetilde P}_n(x,x)$ is strictly
positive for all $x$, except for $n$ odd and $x = 0$. As a
consequence, we see that \eqref{Pk} implies that
$$
 \mathbf{P}_n (x,x)  - \frac{\lambda}{\lambda_n}\, [\mathbf{K}_n(x,c)]^2 +
\frac{\lambda}{\lambda_{n-1}}\, [\mathbf{K}_{n-1}(x,c)]^2 \ge 0
$$
for all $x \in \mathbb{R}^d$ and for all $\lambda > 0$. Recall that $\lambda_n =
1+ \lambda \,\mathbf{K}_{n}(c,c)$. Taking the limit $\lambda \to \infty$, we obtain an
inequality which we state as a proposition.

\begin{prop} For $n\ge1$, and $x, c\in \mathbb{R}^d$,
\begin{equation}\label{ine2}
   \mathbf{P}_{n}(x,x) + \frac{[\mathbf{K}_{n-1}(x,c)]^2}{\mathbf{K}_{n-1}(c,c)}
    \ge \frac{ [\mathbf{K}_n(x, c)]^2}{\mathbf{K}_n(c,c)}.
\end{equation}
\end{prop}

If $x = c$, then the two sides of \eqref{ine2} are equal. This
inequality tends out to be non--trivial even in the case of one
variable. Notice that in one variable, $\mathbf{P}_{n}(x,x) =
[p_n(x)]^2$, where $p_n$ is the orthonormal polynomial. Let us
specify the inequality \eqref{ine2} in the cases of classical
orthogonal polynomials of Jacobi and Laguerre. We use the standard
notation $P_n^{(\alpha,\beta)}$, $\alpha,\beta  > -1$,
 for Jacobi polynomials, which are orthogonal with respect to $(1-x)^\alpha(1+x)^\beta$
 on $[-1,1]$, and $L_n^{(\alpha)}$, $\alpha >-1$, for the Laguerre polynomials,
 which are orthogonal with respect to $x^\alpha e^{-x}$ on $[0,\infty)$. Using the
 formulas in \cite{Sz}, especially (4.3.3), (4.5.3) and  (5.1.1), (5.1.13), the inequality
 (\ref{ine2}) for the Jacobi polynomials with $c=1$, and for the Laguerre polynomials
 with $c =0$, respectively,  becomes the following:

\begin{prop}
For $n \ge 1$ and $x \in [-1,1]$,
$$
 \frac{[P^{(\alpha,\beta)}_n(x)]^2}{P^{(\alpha,\beta)}_n(1)}
+ \frac{n+\beta}{2n + \alpha + \beta +1}\,\frac{[P^{(\alpha+1,\beta)}_{n-1}(x)]^2}
     {P^{(\alpha+1,\beta)}_{n-1}(1)}
 \ge \frac{n+\alpha + \beta + 1}{2n + \alpha + \beta +1}\,
   \frac{[P^{(\alpha+1,\beta)}_n(x)]^2}{P^{(\alpha+1,\beta)}_n(1)}.
$$
For $n \ge 1$ and $x \in [0, \infty)$,
$$
 \frac{[L^{(\alpha)}_n(x)]^2}{L^{(\alpha)}_n(0)} + \frac{[L^{(\alpha+1)}_{n-1}(x)]^2}
    {L^{(\alpha+1)}_{n-1}(0)}  \ge \frac{[L^{(\alpha+1)}_n(x)]^2}{L^{(\alpha+1)}_n(0)}.
$$
\end{prop}

As far as we are aware, these inequalities are new. For example, in the case of
Chebyshev polynomials or $\alpha =\beta = -1/2$ in the Jacobi polynomials, the
inequality becomes
$$
   2 \cos^2 n \theta + \frac{1}{2n-1}\left(\frac{\sin (n-\frac 12)\theta}
       {\sin \frac{\theta}{2}}\right)^2
     \ge \frac{1}{2n +1}\left(\frac{\sin (n+\frac 12)\theta} {\sin \frac{\theta}{2}}\right)^2,
   \quad 0 \le \theta \le \pi.
$$

\section{An example: Krall--type orthogonal polynomials in the unit ball}

Let $B^d$ denote the unit ball of $\mathbb{R}^d$. We consider the inner product
$$
 \langle f, g \rangle_\mu = c_{\mu} \int_{B^d} f(x) g(x)  (1-\|x\|^2)^{\mu-1/2}
 dx,
$$
where $\mu > -1/2$, and $c_{\mu} = {\Gamma(\mu + \frac{d+1}{2})}/
({\pi^{d/2} \Gamma(\mu + \frac{1}{2})})$ is the normalization
constant so that $ \langle 1,1 \rangle_\mu =1$. As an example for
our general results, we add the mass point at the origin, and
consider the inner product
\begin{equation} \label{ball-ip}
 \langle f, g \rangle =  \langle f, g \rangle_\mu + \lambda f(0) g(0), \qquad \lambda > 0.
\end{equation}
We now use \eqref{ex-expl} in Theorem \ref{main-thm} to find an orthogonal basis for
$\langle \cdot, \cdot\rangle$.

Let $\mathcal{H}_n^d$ denote the space of spherical harmonic polynomials of
degree $n$ in $d$ variables. Let $Y_{\nu}^n$, $1 \le \nu \le \dim \mathcal{H}_n^d$,
be an orthonormal basis for $\mathcal{H}_n^d$ in the following. An orthonormal basis for
$\langle f, g \rangle_\mu$ is given explicitly by
(\cite[p. 39]{DX})
\begin{equation} \label{ball-base}
 P_{j,\nu}^n(x) = \left[h_{j,\nu}^{n}\right]^{-1}
 p_j^{(\mu-\frac{1}{2},n-2j+\frac{d-2}{2})}(2\|x\|^2-1)Y_{\nu}^{n-2j}(x), \quad 0 \le j \le n/2,
\end{equation}
where $p_j^{(\mu,n-2j+\frac{d-2}{2})}$ denote the orthonormal Jacobi polynomials,
and $h_{j,\nu}^n$ is the normalizing constant given by
$\left[h_{j,\nu}^n\right]^2 =  (d/2)_{n-2j} / (\mu + \frac{d+1}{2})_{n-2j}$, in which
$(a)_k := a(a+1) \ldots (a+ k-1)$ denotes the shifted factorial.
Since $Y_{\nu}^{n-2j}$ is a homogeneous polynomial, $Y_{\nu}^{n-2j}(0) =0$ unless its
degree is zero, that is, unless $n =2j$. Consequently, $P_{j,\nu}^n(0) = 0$ unless
$j = n/2$ and $n$ is even. Notice that when $n = 2j$, $h_{j,\mu}^n  =1$. Hence, it follows
that $P_{\frac{n}{2},\nu}^{n}(0) = p^{(\mu-\frac{1}{2},\frac{d-2}{2})}_{\frac{n}{2}}(-1)$
if $n$ is even and $P_{j,\nu}^n(0) = 0$ in all other cases. Consequently,  if we define polynomials $Q_{j,\mu}^n$ by
\begin{equation} \label{ball-base2}
Q_{j,\nu}^n(x)= \begin{cases}
  P_{\frac{n}{2},\nu}^n(x) -  \rho_n \mathbf{K}_{n-1}(x,0), & \hbox{if $n$ is even} \\
         P_{j,\nu}^n(x),  & \hbox{otherwise} \end{cases},
\end{equation}
where $\rho_n: =  \lambda  p^{(\mu-\frac{1}{2},\frac{d-2}{2})}_{n/2}(-1)/
 (1 + \lambda \mathbf{K}_{n-1}(0,0))$, then according to \eqref{ex-expl} in
Theorem \ref{main-thm},  $\{Q_{j,\mu}^n : 1 \le \nu \le \dim \mathcal{H}_{n-2j}^d,
0\le 2 j \le n\}$ constitutes an orthogonal
basis with respect to the inner product (\ref{ball-ip}).

To make the expression for transparent, we assume $\mu \ge 0$ and make use of
the following explicit formula for the reproducing kernel $\mathbf{K}_{n}(\cdot,\cdot)$
in \cite{X99},
\begin{align*}
  \mathbf{K}_{n}(x,y) = A_n^\mu \int_{-1}^1 P_n^{(\mu+\frac{d}{2}, \mu + \frac{d}{2} -1)}
           (x \cdot y + \sqrt{1-\|x\|^2} \sqrt{1-\|y\|^2} \, t) (1-t^2)^{\mu-1} dt,
\end{align*}
where  $\mu > 0$ (see \cite{X99} for the case $\mu =0$) and
$$
A_n^\mu := \frac{2\Gamma(\mu + \frac12)\Gamma(\mu+\frac{d+2}{2})\Gamma(n+2\mu + d)}
      {\pi^{1/2} \Gamma(\mu) \Gamma(n+\mu+\frac{d}{2}) \Gamma(2\mu + d+1)}.
$$
We set $y = 0$ in this formula and follow through a sequence of manipulations of
formulas. First we use \cite[(4.5.3)]{Sz} to write $P_n^{(\mu+\frac{d}{2},
\mu+ \frac{d}{2} -1)}(z)$ as a sum of Gegenbauer polynomials, so that we can apply
\cite[Theorem 1.5.6]{DX} to get ride of the integral (for the terms of odd degree
Gegenbauer polynomials, the integrals are automatically zero), the result is a sum of
${\lfloor \frac{n}{2}\rfloor}$ terms of Jacobi polynomials upon using the first formula
on \cite[p. 27]{DX},  which we can use \cite[(4.5.3)]{Sz} again to sum up. The final result
is the following identity,
\begin{align} \label{Kn(x,0)}
\mathbf{K}_n(x,0) =
 \frac{(\mu+\frac{d-1}{2})_{\lfloor \frac{n}{2} \rfloor}}
      {(\mu+\frac{1}{2})_{\lfloor \frac{n}{2} \rfloor}}
    P^{(\frac{d}{2},\mu-\frac{1}{2})}_{\lfloor \frac{n}{2}\rfloor}(1-2\|x\|^2),
\end{align}
where $\lfloor x \rfloor$ denotes the integer part of $x$. In particular, setting $x =0$ gives
\begin{align} \label{Kn(0,0)}
\mathbf{K}_{n}(0,0) = \frac{(\mu+\frac{d-1}{2})_{\lfloor \frac{n}{2}\rfloor}}
    {(\mu+\frac{1}{2})_{\lfloor \frac{n}{2}\rfloor} }
  P^{(\frac{d}{2},\mu-\frac{1}{2})}_{\lfloor \frac{n}{2}\rfloor}(1) =
  \frac{(\mu+\frac{d-1}{2})_{\lfloor \frac{n}{2}\rfloor}}{(\mu+\frac{1}{2})_{\lfloor \frac{n}{2}\rfloor}}\binom{\lfloor \frac{n}{2}\rfloor+ \frac{d}2}{\lfloor \frac{n}{2}\rfloor}.
\end{align}
Substituting these formulas into \eqref{ball-base2} gives a basis of explicit orthogonal
polynomials with respect to the inner product in \eqref{ball-ip}.
Furthermore, using \eqref{Kn(x,0)} and  \eqref{Kn(0,0)} in the formula \eqref{kernel},
we obtain a compact formula for the reproducing kernel
$\mathbf{\widetilde K}_{n}(\cdot,\cdot)$ associated with $\langle \cdot,\cdot\rangle$
in \eqref{ball-ip}. We sum up these results as a proposition.

\begin{prop}
For the inner product $\langle \cdot,\cdot\rangle$ in \eqref{ball-ip}, the polynomials
$Q_{j,\nu}^n$, $1 \le \nu \le \dim \mathcal{H}_{n-2j}^d$, $0 \le j \le n$,  in
\eqref{ball-base2} form an orthogonal basis of degree $n$. Furthermore, the reproducing
kernel of $\Pi_n^d$ with respect to $\langle \cdot,\cdot\rangle$ is given by
$$
   \mathbf{\widetilde K}_{n}(x,y) =  \mathbf{K}_{n}(x,y) - d_n
     P^{(\frac{d}{2},\mu-\frac{1}{2})}_{\lfloor \frac{n}{2}\rfloor}(1-2\|x\|^2)
     P^{(\frac{d}{2},\mu-\frac{1}{2})}_{\lfloor \frac{n}{2}\rfloor}(1-2\|y\|^2),
$$
where the constant $d_n$ is given by
$$
  d_n = \frac {\lambda} {1 + \lambda \mathbf{K}_{n}(0,0)}
     \left[\frac{(\mu+\frac{d-1}{2})_{\lfloor \frac{n}{2}\rfloor}}
    {(\mu+\frac{1}{2})_{\lfloor \frac{n}{2}\rfloor} }\right]^2.
$$
\end{prop}

\end{document}